\documentclass[11pt]{article}
\usepackage{url}
\begin{document}
\bibliographystyle{plain}
\baselineskip=15pt
\title{What is the Inverse of Repeated Square and Multiply Algorithm?}
\author{H. Gopalkrishna Gadiyar$^{(1)}$, K M Sangeeta Maini$^{(1)}$,\\
R. Padma$^{(1)}$ and Mario Romsy$^{(2)}$\\
~~\\
$^{(1)}$AU-KBC Research Centre,\\ M. I. T. Campus of Anna University,\\  Chromepet, Chennai 600 044, INDIA\\
E-mail: \{gadiyar, padma\}@au-kbc.org\\
kmsmlm\_2001@rediffmail.com \\
~~\\
$^{(2)}$ Universit\"{a}t der Bundeswehr München,\\
Institut f\"{u}r theoretische Informatik und \\
Mathematik, 85577 Neubiberg, Germany\\ 
E-mail: mario.romsy@unibw.de}
\date{~~~}
\maketitle

\begin{abstract}
It is well known that the repeated square and multiply algorithm is an efficient way of modular exponentiation. The obvious question to ask is if this algorithm has an inverse which would calculate the discrete logarithm and what is its time compexity. The technical hitch is in fixing the right sign of the square root and this is the heart of the discrete logarithm problem over finite fields of characteristic not equal to 2. In this paper a couple of probabilistic algorithms to compute the discrete logarithm over finite fields and their time complexity are given by bypassing this difficulty. One of the algorithms was inspired by the famous $3x+1$ problem.  

\noindent {\bf Key words.} Discrete logarithm, Legendre symbol, 3x+1 problem.

\end{abstract}

\section{Introduction}
Let $p$ be an odd prime number and let $a$ be a primitive root. Then any $b \in (Z/pZ)^*$ can be expressed as 
\begin{equation}
b \equiv a^n {\rm ~mod~}p \label{eq:b}
\end{equation}
for a unique integer $n$ with $1 \le n \le p-1$. $n$ is called the index or discrete logarithm of $b$ to the base $a$ modulo $p$. The discrete logarithm problem over prime fields is to find $n$ given $a$ and $b$ modulo $p$. When $p$ is sufficiently large and random, the discrete logarithm problem is believed to be computationally difficult and hence is the basis of the security of many cryptographic algorithms like the Diffie - Hellman key exchange protocol, El Gamal Cryptosystem, El Gamal signature scheme etc. The well known algorithms to compute the discrete logarithm problem are the baby step - giant step method, Pollard's rho method, Pohlig - Hellman method and the index calculus method \cite{menezes:oorschot:vanstone}, \cite{mccurley}, \cite{odlyzko}, \cite{schirokauer:weber:denny}. In this the first two algorithms are square root algorithms and the Pohlig - Hellman method works efficiently for those primes $p$ for which all the prime factors of $p-1$ are small. The index calculus method is a subexponential time algorithm. All these algorithms find $n$ modulo $p-1$. There is an analysis of the discrete logarithm problem using $p$-adic methods in \cite{gadiyar:maini:padma}. In \cite{riesel}, a straightforward formula for solving the discrete logarithm problem in certain cases is obtained using the Fermat quotient and its generalizations. For an extensive bibliography on this problem look at the website \url{http://www.cs.uwaterloo.ca/~shallit/bib/dlog.bib}. The Pohlig-Hellman reduction is not itself a method to solve the discrete logarithm problem, it reduces the problem to a number of discrete logarithm problems in groups of prime order. So in practice one first applies the Pohlig-Hellman reduction and then for example the Pollard-$\rho $- method.So far, there have been no known polynomial time algorithm to compute the discrete logarithm for a random prime $p$. 

We know that the modular exponentiation $a^n {\rm ~mod~}p~(=b)$ in (\ref{eq:b}) is performed efficiently using the repeated square and multiply algorithm. In this paper we ask the question of what is the inverse of this algorithm and how much time does it take. As an answer to this question we give a couple of probabilistic algorithms that compute the discrete logarithm. The Legendre symbol of $b$ determines the least significant bit (l.s.b.) of $n$ which is the index of $b$ to the base $a$ in (\ref{eq:b}). When 
the l.s.b. of $n$ is $0$, the next least significant bit of $n$ is 
obtained by extracting the `correct' square root. The last $r$ 
significant bits can be unambiguously and efficiently determined, if $p-1 =2^rs$, where $s$ is odd and $r \ge 1$ \cite{peralta}. The ambiguity 
starts from the last $r+1^{th}$ significant bit onwards. 

In this paper we give two probabilistic algorithms which bypass this 
ambiguity and compute the discrete logarithm over prime fields. The first algorithm can be thought of as a randomized inverse of the repeated 
square and multiply algorithm. The second algorithm was inspired by the 
famous $3x+1$ problem. 
The algorithms are immediately extendable to other  finite fields (including finite fields of characteristic 2 with a slight modification). In Section \ref{sect:legendre} we explain how the properties of the Legendre symbol can be used to determine the last $r$ significant bits of the index $n$ and also give a time estimate for computing square roots modulo $p$. In Section \ref{sect:inversesqmul} we state our main algorithm and give numerical examples over a prime field. A modification of the algorithm for finite fields of characteristic 2 is also given along with a couple of examples. This algorithm is extendable to the elliptic curve discrete logarithm problem. In Section \ref{sect:3x} we state the $3x+1$ problem and a variant of the algorithm given in Section \ref{sect:inversesqmul} is presented with an example. In Section \ref{sect:timecomplexity} we analyze the time complexity of the algorithms.

\section{The Legendre Symbol and Square Roots}
\label{sect:legendre}

The Legendre symbol $\left(\frac{x}{p}\right )$, for any integer $x$, 
with $(x,p)=1$ is defined as follows.
\begin{eqnarray}
\left (\frac{x}{p}\right )=\left\{ \begin{array}{rl}
1 \, , &\mbox{if $x^{\frac{p-1}{2}} \equiv 1$~mod~$p$}\\
-1 \, , &\mbox{if $x^{\frac{p-1}{2}} \equiv -1$~mod~$p$}\, . 
\end{array}
\right. \label{eq:legx}
\end{eqnarray}
The definition (\ref{eq:legx}) of Legendre symbol can be restated as follows. $\left (\frac{x}{p}\right )=1$, if $x$ is a quadratic residue (that is, a square) modulo $p$ and is equal to $-1$, if $x$ is a quadratic non-residue (that ~is, a non-square) modulo $p$. Since $a$ is a primitive root, 
\begin{equation}
\left (\frac{a}{p}\right )\equiv a^{\frac{p-1}{2}}\equiv -1 {\rm~mod~} p \, . \label{eq:lega}
\end{equation}
Using the property of the Legendre symbol
\begin{equation}
\left (\frac{xy}{p}\right ) =\left (\frac{x}{p}\right )\left (\frac{y}{p}\right ) \, ,
\end{equation}
one has 
\begin{equation}
\left (\frac{b}{p}\right )= \left (\frac{a^n}{p}\right )=\left (\frac{a}{p}\right )^n=(-1)^n.
\end{equation}
Thus the Legendre symbol of $b$ determines whether $n$ is odd or even. In other words, the Legendre symbol of $b$ determines the least significant bit of $n$. Let us write $n$ in its binary representation
\begin{equation}
n=n_0+2n_1+2^2n_2+\cdots +2^kn_k \, , \label{eq:bin}
\end{equation}
where each $n_i$ is 0 or 1. If $n$ is even (or odd), then $n_0=0~ ({\rm or~} 1)$. The next bit $n_1$ is determined by dividing $n~({\rm or}~n-1) $ by 2 and checking whether $\frac{n}{2}$ (${\rm or}~\frac{n-1}{2}$ ) is even or odd. In terms of $a$ and $b$ modulo $p$, this amounts to finding the `correct' square root of $b$ (${\rm or}~\frac{b}{a}$) modulo $p$.

Let us assume that $b$ is a quadratic residue modulo $p$. Then the square roots of $b$ are $b^{\frac{1}{2}}$ and $-b^{\frac{1}{2}}$ and hence from (\ref{eq:lega}), the index of the square roots to the base $a$ are $\frac{n}{2}$ and $\frac{p-1}{2}+\frac{n}{2}$ modulo $p-1$. If $p\equiv 1 {~\rm mod}~4$, then the l.s.b. of $\frac{p-1}{2}$ is 0. Hence the l.s.b. of $\frac{n}{2}$= the l.s.b. of  $\frac{p-1}{2}+\frac{n}{2}$. In other words, the Legendre symbol of $b^{\frac{1}{2}}$ or $-b^{\frac{1}{2}}$ will determine the value of the bit $n_1$ unambiguously. More generally, if $p-1=2^rs$, where $s$ is odd and $r\ge 1$, then the last r significant bits of $n$ can be unambiguously determined. See for example \cite{peralta}.
The difficulty arises from the $(r+1)^{th}$ least significant bit onwards.

Let us explain this with the case of $r=1$. In this case $p\equiv 3 \rm{~mod}~4$. Hence the l.s.b. of $\frac{p-1}{2}$ is 1. Now,  if the l.s.b. of  $\frac{n}{2}=1$, then the l.s.b. of  $\frac{p-1}{2}+\frac{n}{2}=0$ and vice versa. Thus it is not possible to determine the correct value of the bit $n_1$. In the next section we give our main algorithm which bypasses this difficulty.

Note that the Legendre symbol can be calculated in polynomial time ($O(\log ^2p)$) \cite{menezes:oorschot:vanstone}. If $p \equiv ~3 \rm{~mod}~4$, then $r=1$, and in this case the square roots of a quadratic residue $x$ are $x^{\frac{s+1}{2}} ~{\rm mod} ~p$ and $-x^{\frac{s+1}{2}} ~\rm{mod} ~p$. Hence the square root can be calculated in polynomial time ($O(\log ^3p)$). If $p \equiv ~1 {\rm ~mod}~4$, then there is a polynomial time algorithm ($O(\log ^4p)$) to compute a square root of a quadratic residue, provided we could find a quadratic nonresidue modulo $p$ (\cite{menezes:oorschot:vanstone}, \cite{bach:shallit}). Since $a$ is a primitive root, it is a quadratic nonresidue and hence the square roots of a quadratic residue $x$ modulo $p$ can be calculated deterministically in polynomial time. 

\section{The Repeated Square and Multiply Algorithm and Its Inverse}
\label{sect:inversesqmul}

In this section we give our main algorithm. In Section \ref{subsect:modexp} we explain how modular exponentiation is done using the repeated square and multiply algorithm and describe the difficulty in trying to invert the process. In Section \ref{subsect:alg} we give our main algorithm for computing the discrete logarithm.  Section \ref{subsect:exp} explains the algorithm and Section \ref{subsect:eg} gives two examples. In Section \ref{subsect:finite} we generalize the algorithm to all finite fields including finite fields of characteristic 2.  This section also contains examples of computing discrete logarithm over finite fields of characteristic 2. Finally, Section \ref{subsect:elliptic} contains a note on the elliptic curve discrete logarithm problem over binary fields.

\subsection{Modular Exponentiation and its Inverse}
\label{subsect:modexp}

The repeated square and multiply algorithm is used to compute modular exponentiation in polynomial time. Let us quickly recall how we compute $b$ given $a$, $n$ and $p$ as in (1) using this algorithm. If $n$ has the binary representation as in (\ref{eq:bin}), then let $a_0=a$, $b_0=1$ and inductively compute $a_j=a_{j-1}^2$ mod $p$ and $b_j=b_{j-1}.a_j$ mod $p$, if $n_j=1$ and $b_j=b_{j-1}$, if $n_j=0$, for $j=1$ to $k$. Then $b_k$ is the value of $a^n \rm{~mod~}p$. That is,
\begin{equation}
b \equiv  a^n \equiv \left (a^{2^k} \rm{~mod~}p \right )^{n_k} \cdots  \left (a^{2^1} \rm{~mod~}p \right )^{n_1}~ \left (a^{2^0} \rm{~mod~}p \right )^{n_0} \, {\rm ~mod~}p \, .
\end{equation} 
It is clear that the inverse of this algorithm is to divide and repeatedly extract square root. Division is done when the bit $n_i$ is $1$, just as multiplication is done in the repeated square and multiply algorithm when $n_i$ is 1. Note that knowing the bits $n_i$'s is equivalent to knowing the value of $n$. Also, if the `correct' square root can be taken every time, it will fix the correct value of the bit $n_i$ and hence $n$ can be calculated in polynomial time. But the difficulty is in fixing the correct square root. We know from Section \ref{sect:legendre} that the last $r$ significant bits $n_i$ can be unambiguously determined in polynomial time by dividing by $a$ if the Legendre symbol is $-1$ and extract any one square root or by just extracting any one square root if the Legendre symbol is $1$. From the last $r+1^{th}$ significant bit onwards, we do not know which is the right bit.  Now we give our probabilistic algorithm which bypasses this problem. 

\subsection{Algorithm to Compute Discrete Logarithm}
\label{subsect:alg}

\noindent {\bf Input:} $a,~b,~p$, where $a^n \equiv b {\rm ~mod~}p$ \, .\\

\noindent {\bf Output:} $n$\\

\noindent {\bf Step 1.} Choose an integer $B$ and create Table I consisting of ($a^{k_j}$ mod $p$, $k_j$) where $j \le B$. Here 
$\{ k_j \}$ is any subsequence of integers. For example, $k_j=j$ or $k_j=2^j$.\\

\noindent {\bf Step 2.}  Initialize $i\leftarrow 1$, $k\leftarrow 0$, $l\leftarrow 1$ and $m\leftarrow n$, $b[1] \leftarrow b$, $m_1[1] \leftarrow n$, $c_1[1] \leftarrow b$, $c_2[1] \leftarrow b$ and $m_2[1] \leftarrow n$. 
\begin{tabbing}
{\bf Step 3.}  \= ~i) If $\left (\frac{b}{p}\right )=-1$ then goto Step 4.\\
\> ii) If $\left (\frac{b}{p}\right )= 1$ then goto Step 6.
\end{tabbing}

\begin{tabbing}
{\bf Step 4.}  \= ~~i) $b\leftarrow \frac{b}{a} {\rm ~mod~}p$ and $m\leftarrow m-1$. \\
\> ~ii) Goto Step 5. \\
\> iii) If Step 5 does not solve for $n$, $i\leftarrow i+1$.\\
\> ~~~~~~~store $b[i]\leftarrow b$ and $m_1[i]\leftarrow m$ in Table II. \\
\> ~~iv) Goto Step 6.
\end{tabbing}

\begin{tabbing}
{\bf Step 5.} \= ~~i) If $b \equiv a^{k_j}{\rm ~mod~}p$ for any $j \le B$ in Table I, Solve($m,k_j,k$).\\

\> ~ii) If  $b\equiv b[j]{\rm ~mod~}p$ for any $j$ in Table II,  Solve($m,m_1[j],k$).\\

\> iii)  If $b \equiv c_1[j] {\rm ~or~} c_2[j]{\rm ~mod~}p$ for any $j$ in Table III, Solve($m, m_2[j],k$).
\end{tabbing}

\begin{tabbing}
{\bf Step 6.} \= ~~i) $b\leftarrow b^{\frac{1}{2}} {\rm ~mod~} p$ and $m\leftarrow \frac{m}{2}$. $k\leftarrow k+1$. Goto Step 5. \\
\> ~ii) If Step 5 does not solve for $n$, $b\leftarrow p-b {\rm ~mod~} p.$ Goto Step 5. \\
\> iii) If Step 5 does not solve for $n$, $l\leftarrow l+1$, \\
\> ~~~~~~store $c_1[l]\leftarrow b$, $c_2[l]\leftarrow p-b$ and $m_2[l]\leftarrow m$ in Table III. \\
\> ~iv) $b \leftarrow c_1[l] {\rm ~or~} c_2[l]$ randomly.\\
 \> ~~v) Goto Step 3.
\end{tabbing}

\noindent {\bf Solve()}

\noindent Solve($x,y,t$): Solve the linear congruence: 
\begin{equation}
2^t x \equiv 2^t y {\rm ~mod~}p-1 \, . 
\end{equation}

\noindent Return $n$

\subsection{Explanation of the Algorithm}
\label{subsect:exp}

Table I consists of $B$ precomputed powers of $a$ and the corresponding discrete logarithms. Table II consists of intermediate values of $b$ and the corresponding values of $m$ as a function of $n$ when the Legendre symbol is -1. Table III consists of the intermediate values of $b$ which are square roots and the corresponding values of $m$ as a function of $n$ when the Legendre symbol is 1. Since $m$ is replaced by $m-1$ or $\frac{m}{2}$, $m$ is always a linear function of $n$. The algorithm is probabilistic as we select one of the square roots randomly in Step 6.

The final step is to solve a linear congruence modulo $p-1$, if the new value of $b$ matches any of the integers in Table I, II or III. If $b$ matches a value in Table I, then $x$ will be a linear function of $n$ and $y$ will be a constant. In other cases, both $x$ and $y$ will be linear functions of $n$. Note that a linear congruence can be solved in polynomial time.

When $b$ coincides with a value in Table I, the value of $n$ can be uniquely obtained by solving the linear congruence.

When $b$ coincides with a value in Table II or III, then the corresponding value of $n$ can be got modulo ($\frac{p-1}{d}$), for a divisor $d$ of $p-1$. Hence there will be $d$ solutions modulo $p-1$ and we have to choose the correct value of $n$ modulo $p-1$ by trial and error. If $d$ is too big, then one can start the algorithm again from somewhere in the middle of the tree where we can choose the other square root.   

$k$ counts the number of times we take square roots modulo $p$. Note that while solving the linear congruence, we multiply both sides by $2^k$, so that the denominator of $m$ gets cleared (as $2$ is not invertible modulo $p-1$.) This also takes care of the fact that in Table III, though we store two square roots, the exponent $m$ is taken to be $\frac{m}{2}$, as whether we take $\frac{m}{2}$ or $\frac{m}{2}+\frac{p-1}{2}$, in Solve(), multiplication by $2^k$ would remove this ambiguity.

\subsection{Examples}
\label{subsect:eg}

In this section we explain the algorithm in Section \ref{subsect:alg} with a small prime. Let $p=103$. $a=5$ is a primitive root of $p$.

\noindent{\bf Example 1.}  This is an example of collision with an element in Table I. Let $b=84$. Let $B=7$ and $k_j=2^{j-1} {\rm ~mod~}p$ for $j=1,\cdots 7$. 

\begin{center}
{\bf Table I}\\

\vspace{0.25cm}
\begin{tabular}{|c|c|c|c|c|c|c|c|}
 \hline
  $~j~$& $~0~$ & $~1~ $ & $~2~$ &~3 ~  & $~4~$&~$5$~ & ~$6$~\\
 \hline
 $5^{2^j}{\rm ~mod~}103$&5 &25  & 7&49 &32&97&36 \\
\hline
 \end{tabular}
\end{center}

\begin{center}
{\bf Discrete Logarithm Calculation for $b=84$}\\

\vspace{0.25cm}
\begin{tabular}{|c|c|c|c|c|c|}
 \hline
  $~b~$& $~\left (\frac{b}{103}\right )~$ & $~\frac{b}{a}{\rm ~mod~}103~ $ & $~b^{\frac{1}{2}},~ -b^{\frac{1}{2}} {\rm ~mod~}103$ &~random ~  & $~m~$ \\
~&~&~&~&sqrt&~\\
 \hline
 84&-1 &58  & -- &-- &$n-1$ \\
58&1 &--  &26,77  & 77&$\frac{n-1}{2}$ \\
77&-1 &36  & --&-- &$\frac{n-1}{2}-1$ \\
\hline
 \end{tabular}
\end{center}
Since $36\equiv 5^{2^6} {\rm ~mod~} 103$ and $k=1$ as we have taken square root only once, after multiplying both sides by $2$ we get the congruence
\begin{equation}
n-3\equiv 2^7 {\rm ~mod~}102
\end{equation}
and thus $n \equiv 29 {\rm ~mod~}102$.

Note that the binary digits of 29 are given by (11101). Comparing these bits and the second column of the above table, we find that we have wrongly chosen the square root 77 in the second row, yet we are lucky to find the discrete logarithm at the third step itself, while if we had taken the correct path, it would have taken us 6 steps. There is a trading off between extraction of square roots and the precomputation of powers of $a$.
 
Note that Table II in the algorithm corresponds to third and sixth columns of the above table and Table III corresponds to fourth and sixth columns.

\noindent{\bf Example 2.} This example gives a collision in Table II or III.

\begin{center}
{\bf Discrete Logarithm Calculation for $b=99$}\\
\vspace{0.25cm}
\begin{tabular}{|c|c|c|c|c|c|}
 \hline
  $~b~$& $~\left (\frac{b}{103}\right )~$ & $~\frac{b}{a}{\rm ~mod~}103~ $ & $~b^{\frac{1}{2}},~ -b^{\frac{1}{2}} {\rm ~mod~}103$ &~random ~  & $~m~$ \\
~&~&~&~&sqrt&~\\
 \hline
 99&-1 &61  & -- &-- &$n-1$ \\
61&1 &--  &24,79  & 24&$\frac{1}{2}(n-1)$ \\
24&-1 &46  & -- &-- &$\frac{1}{2}(n-3)$ \\
46&1 &--  &47,56  & 56&$\frac{1}{4}\left (n-3\right )$ \\
56&1&--&46,57&--&$\frac{1}{8}\left (n-3\right )$\\ 
\hline
 \end{tabular}
\end{center}
Note that $46$ in the fourth column matches the $46$ in the third column. Equating the corresponding values of $m$, multiplying both sides by $2^3$ gives
\begin{equation}
n-3 \equiv 4(n-3) {\rm ~mod~} 102 \, .
\end{equation}
Solving the linear congruence gives $n \equiv 3 {\rm ~mod~}34$. Hence there are $3$ possible values for $n$ mod 102, namely, 3, 37 and 71. One can easily check that $37$ is the correct value of $n$.  

\subsection{Discrete Logarithm over Finite Fields}
\label{subsect:finite}

It is clear that the algorithm given in Section \ref{subsect:alg} is, just as it is, extendable to finite fields of characteristic $p > 3$, as the analogue of Legendre symbol and efficient computation of square roots exist in these fields \cite{bach:shallit}. 

When the characteristic of the finite field is $2$, every element in the field is a square and every element has exactly one square root. \cite{fong:hank:lopez:menezes} gives an efficient algorithm for computing square roots over finite fields of characteristic 2. Note that if we choose a normal basis, then the square root operation reduces to a mere cyclic shift \cite{menezes:oorschot:vanstone}. Hence, our algorithm in Section \ref{subsect:alg} can be modified in this case as follows. Step 3 should randomly decide whether the l.s.b. of $m$ is 1 or 0. That is,

\begin{tabbing}
{\bf Step 3.}  \= ~0) Randomly choose a bit 0 or 1.\\
\> ~i) If the bit is 1 then goto Step 4.\\
\> ii) If the bit is 0 then goto Step 6.
\end{tabbing}

In Step 6, as there is only one square root, we need not perform ii) and Table III will consist of only $c_1[l]$ and the corresponding value of $m$. Also Step iv) should be skipped. 

To be precise, the randomness in our algorithm for $p>3$ in the selection of the square root in iv) of Step 6 has been shifted to Step 3 where we randomly fix the l.s.b. of $m$ for characteristic 2 fields. 

We give two toy examples below.

Let us consider the finite field $F_{2^7}$ with the primitive polynomial $f(x)=x^7+x+1$. Thus, $x$ is the generator of the multiplicative group $F_{2^7}^{*}$ of the finite field $F_{2^7}$.

\noindent{\bf Example 1.} Let us create Table I with $B=7$ and $k_j=2^j$ for $j=0 \cdots 6$. Let $b=x^4+x^3+x^2+1$.

\begin{center}
{\bf Table I}\\

\vspace{0.25cm}
\begin{tabular}{|c|c|c|c|c|c|c|c|}
 \hline
  $~j~$& $~0~$ & $~1~ $ & $~2~$ &~3 ~  & $~4~$&~$5$~ & ~$6$~\\
 \hline
 $x^{2^j}{\rm ~mod~}f(x)$&$x$ &$x^2$ & $x^4$&$x(x+1)$&$x^2(x^2+1)$ &$x(x^3+x+1)$&$x(x^3+1)$ \\
\hline
 \end{tabular}
\end{center}

\begin{center}
{\bf Discrete Logarithm Calculation for $b=x^4+x^3+x^2+1$}\\

\vspace{0.25cm}
\begin{tabular}{|c|c|c|c|c|}
 \hline
  $~b~$& ~random ~ & $~\frac{b}{x}{\rm ~mod~}f(x)~ $ & $~b^{\frac{1}{2}}{\rm ~mod~}f(x)$  & $~m~$ \\
~&bit&~&~&~\\
 \hline
 $x^4+x^3+x^2+1$&0 &--  & $x^5+x+1$ &$\frac{n}{2}$ \\
$x^5+x+1$&1 &$x^4(x^2+1)$ &-- &$\frac{n}{2}-1$ \\
$x^4(x^2+1)$ &1&$x^3(x^2+1)$ &--&$\frac{n}{2}-2$ \\
$x^3(x^2+1)$ &1&$x^2(x^2+1)$ &--&$\frac{n}{2}-3$ \\
\hline
 \end{tabular}
\end{center}
From Table I, $x^2(x^2+1)\equiv x^{2^4}{\rm ~mod~}f(x)$. Hence,
\begin{equation}
\frac{n}{2}-3\equiv 2^4 {\rm ~mod~}127 \, ,
\end{equation}
which gives $n \equiv 38 {\rm ~mod~}127$.

\noindent{\bf Example 2.} Let $b=x^6+x^5+x^3+x+1$.
\begin{center}
{\bf Discrete Logarithm Calculation for $b=x^6+x^5+x^3+x+1$}\\

\vspace{0.25cm}
\begin{tabular}{|c|c|c|c|c|}
\hline
$~b~$& ~random ~ & $~\frac{b}{x}{\rm ~mod~}f(x)~ $ & $~b^{\frac{1}{2}}
{\rm ~mod~}f(x)$  & $~m~$ \\
~&bit&~&~&~\\
\hline
 $x^6+x^5+x^3+x+1$&0 &--  & $x^6+x^5+x^4+x^2+x+1$ &$\frac{n}{2}$ \\
$x^6+x^5+x^4+x^2+x+1$&1 &$x^6+x^5+x^4+x^3+x$& -- &$\frac{n}{2}-1$ \\
$x^6+x^5+x^4+x^3+x$&1&$x^5+x^4+x^3+x^2+1$&--&$\frac{n}{2}-2$ \\
$x^5+x^4+x^3+x^2+1$&0&--&$x^6+x^5+x^3+x+1$&$\frac{1}{2}\left (\frac{n}
{2}-2\right )$ \\
\hline
\end{tabular}
\end{center}
The value of the original $b$ and the square root in the fourth row are 
the same. Equating the corresponding values of $m$, we get
\begin{equation}
n \equiv \frac{1}{2}\left (\frac{n}{2}-2\right ) {\rm ~mod~}127 \,.
\end{equation}
Solving this linear congruence gives $n \equiv 41 {\rm ~mod~}127$.

\subsection{Elliptic Curve Discrete Logarithm over Binary Fields}
\label{subsect:elliptic}

A natural question to ask now is if our algorithm can be extended to elliptic curves as well. Since there is no analogue of the Legendre symbol, it is not possible to extend this algorithm to elliptic curves. Note that in any group of odd order (written multiplicatively) every element has exactly one square root. So if we can calculate this square root efficiently, the modified algorithm of this section can be applied. When the elliptic curve is defined over finite fields of characteric 2, there is an efficient algorithm for point-halving \cite{knud}, \cite{schro} if the cardinality of the curve is odd. Hence in this case we can generalize the algorithm given in Section \ref{subsect:finite} to compute the discrete logarithm on such elliptic curves.    

\section{The $3x+1$ Problem and the Discrete Logarithm Problem}
\label{sect:3x}

In this section we give a variant of the main algorithm given in Section \ref{subsect:alg}. This was inspired by the famous $3x+1$ problem. 
The $3x+1$ problem which was posed by L. Collatz, states that if
\begin{eqnarray}
T(x)=\left\{ \begin{array}{ll}
3x+1 \, , &\mbox{if $x \equiv 1$~mod~2}\\
\frac{x}{2}\, , &\mbox{if $x \equiv 0$~mod~2} \, , 
\end{array}
\right.
\end{eqnarray}
then for any positive integer $x$ there exists an integer $k > 0$ with
$T^k(x) = 1$. This problem remains open since 1937. For an annotated bibliography of this problem, see \cite{lagarias}.  

Note that in this problem if $x$ is odd, the function $T$ converts it into an even integer by multiplying $x$ by $3$ and then adding 1, while if $x$ is even, it divides $x$ by $2$. The iteration will terminate once $T^k(x)=2^l$ for some integers $k$ and $l$.

In the algorithm we gave in Section \ref{subsect:alg}, if the Legendre symbol is $-1$ (that is, the index of $b$ is odd), we divided $b$ by $a$ so that the index of the new value of $b$ becomes even and if the Legendre symbol is $1$, (that is, the index of $b$ is even), we calculated the square roots of $b$ so that the index is halved. 

Now it is clear how we are going to modify the algorithm in Section \ref{subsect:alg}. We will assume for the sake of simplicity that $(3,p-1)=1$. If the Legendre symbol is $-1$, then compute $b^3a {\rm ~mod~}p$. That is, in Step 4 (i), we do 
\begin{equation}
b\leftarrow b^3a {\rm ~mod~}p {\rm ~and~} m\leftarrow 3m+1
\end{equation}
and the rest of the algorithm goes as before. 

\subsection{Example} 
\label{subsect:eg1}
We explain our algorithm with an example, again with a small prime. Let us take $p=101$. $a=2$ is a primitive root of 101. Let $b=72$. Let $B=7$ and $k_j=2^j$, for $j=0 \cdots 6$. 

\begin{center}
{\bf Table I}\\

\vspace{0.25cm}
\begin{tabular}{|c|c|c|c|c|c|c|c|}
 \hline
  $~j~$ &$~0~$& $~1~ $ & $~2~$ &~3 ~  & $~4~$&~$5$~ & ~$6$~\\
 \hline
 $2^{2^j} {\rm ~mod~}101$&2 &4  & 16&54&88&68&79 \\
\hline
 \end{tabular}
\end{center}

\begin{center}
{\bf Discrete Logarithm Calculation for $b=72$}\\

\vspace{0.25cm}
\begin{tabular}{|c|c|c|c|c|c|}
 \hline
  $~b~$& $~\left (\frac{b}{101}\right )~$ & $~b^3a{\rm ~mod~}101~ $ & $~b^{\frac{1}{2}},~ -b^{\frac{1}{2}} {\rm ~mod~}101$ &~random ~  & $~m~$ \\
~&~&~&~&sqrt&~\\
 \hline
 72&-1 &5 & --&--&$3n+1$ \\
5&1 &--  &45,56  & 56&$\frac{3n+1}{2}$ \\
56&1 &--  & 37,64 &37&$\frac{3n+1}{4}$ \\
37&1&--&21,80&80&$\frac{3n+1}{8}$ \\
80&1&--&22,79&--&$\frac{3n+1}{16}$ \\
\hline
 \end{tabular}
\end{center}
Since $79  \equiv 2^{2^6} {\rm ~mod~} 101$, and $k=4$, we have
\begin{equation}
3n+1 \equiv 1024 \equiv 24 {\rm ~mod~}100 \, .
\end{equation}
The solution of this linear congruence is given by $n\equiv 41 {\rm ~mod~}100$.

\section{Time Complexity of the Algorithms}
\label{sect:timecomplexity}

Let us first look at the algorithm in Section \ref{subsect:finite} over finite fields of characteristic 2. Here we randomly choose to do division or extract square root. Let $\rho :~ {\cal N}~ \rightarrow ~\{0,1\}$ be the random decision function. Starting with $r_0~=b=a^n$, we have for $i \in \cal{N}$
\begin{eqnarray}
r_i=\left\{ \begin{array}{ll}
\sqrt{r_{i-1}}, &\mbox{~if~ $\rho (i)=0$}\\
\frac{r_{i-1}}{a} ,&\mbox{~if~ $\rho (i)=1$} \, . 
\end{array}
\right.
\end{eqnarray}
Let Table I consist of $B \in \cal{N}$ precomputed values and label them $r_{-B},...,r_{-1}$. Since we need to keep track of the exponents we store pairs of the form $(a^k,~k)$, hence write $(r_{-B},k_{-B}),...,(r_{-1},k_{-1})$, where $\{kj\}$ is the chosen subsequence of integers (as in Section
3.2). So starting with $(r_0,m_0)$, $m_0=n$, we calculate a random walk by
\begin{eqnarray}
(r_i, m_i)=\left\{ \begin{array}{ll}
(\sqrt{r_{i-1}}, \frac{m_{i-1}}{2}), &\mbox{~if~ $\rho (i)=0$}\\
(\frac{r_{i-1}}{a}, m_{i-1}-1) ,&\mbox{~if~ $\rho (i)=1$} \, . 
\end{array}
\right.
\end{eqnarray}
Then we look for a collision $r_i=r_j$ for $i \neq j$. Since $n$ is unknown, the $m_i$ are linear functions in $n$.

Note that in the algorithm given in Section \ref{subsect:alg}, the random function will decide which square root will be taken. Since we store both square roots, we add a control bit $b_i$ in each step where $b_i=1$ means a square root was taken ($b_i=0$ means a division.) In step $i$, when we look for collisions of $r_i$ with some previous element, we test $r_i=r_j$ if $b_j=0$ and $r_i= \pm r_j$ if $b_j=1$. (Note that if $r_j$ is one square root then $-r_j$ is the other and that the associated linear functions are the same.) Hence in this case, we have
\begin{eqnarray}
(r_i, m_i, b_i)=\left\{ \begin{array}{ll}
(\min(\sqrt{r_{i-1}},p-\sqrt{r_{i-1}}),~ \frac{m_{i-1}}{2},1), &\mbox{~if~
$(\frac{r_{i-1}}{p})=1$~and~ $\rho (i)=0$}\\
(\max(\sqrt{r_{i-1}},p-\sqrt{r_{i-1}}),~ \frac{m_{i-1}}{2}, 1), &\mbox{~if~
$(\frac{r_{i-1}}{p})=1$~and~ $\rho (i)=1$}\\
(\frac{r_{i-1}}{a}, m_{i-1}-1, 0) ,&\mbox{~if~ $(\frac{r_{i-1}}{p})=-1$} \, . 
\end{array}
\right.
\end{eqnarray}
With respect to the algorithm given in Section \ref{sect:3x}, only the iteration function is changed. Since this is a random walk, the expected number of steps should be about $O(\sqrt{p})$. In comparison with the Pollard $\rho $ -method, our algorithms are probably slower, since taking square roots and division are (in general) slower than squaring and multiplication. Also, Floyd's cycle detection method will not apply here as the function $\rho $ is a random decision function. On the other hand, $\rho $ being a `random
decision function' gives the walk a random pattern.

\section{Conclusion and Future Directions}
\label{sect:conc}
In this paper we have asked the question of what is the inverse of repeated square and multiply algorithm and given a couple of probabilistic algorithms to compute the discrete logarithm. The algorithms are parallelizable. It is noted that the algorithm given for binary finite fields can also be extended to elliptic curves over such fields. Analysis of the algorithms shows that these algorithms are of square root type. like the baby step-giant step method, Pollard's $\rho $ method etc. \cite{teske1} Though the algorithm given in Section \ref{sect:3x} does not `carry the idea of the inverse of square and multiply', it will be worthwhile to study the more general scheme $x \rightarrow ax+b$, with $(a, p-1)=1$, if $x$ is odd, and $x \rightarrow \frac{x}{2}$, if $x$ is even and see if some interesting algebra could be uncovered.

\end{document}